\documentclass[amsfonts]{article}
\usepackage{amssymb}

\newtheorem{theorem}{Theorem}[section]
\newtheorem{corollary}[theorem]{Corollary}

\newtheorem{proposition}[theorem]{Proposition}
\newtheorem{lemma}[theorem]{Lemma}
\newtheorem{definition}[theorem]{Definition}
\newtheorem{remark}[theorem]{Remark}

\def\bR{\mathbb{R}}
\def\bZ{\mathbb{Z}}

\def\bN{\mathbb{N}}

\def\cB{\mathcal{B}}

\def\cF{\mathcal{F}}

\def\cJ{\mathcal{J}}
\def\cK{\mathcal{K}}

\def\cM{\mathcal{M}}

\begin{document}

\title{Explicit Conditions for the Convergence of \\ Point Processes
Associated to Stationary Arrays}

\author{Raluca Balan\footnote{Corresponding author. Department of Mathematics and Statistics, University of Ottawa,
585 King Edward Avenue, Ottawa, ON, K1N 6N5, Canada. E-mail address:
rbalan@uottawa.ca} \ \footnote{Supported by a grant from the Natural
Sciences and Engineering Research Council of Canada.} \ and Sana
Louhichi\footnote{Laboratoire de math\'ematiques, \'Equipe de
Probabilit\'es, Statistique et mod\'elisation, Universit\'e de
Paris-Sud, B\^at. 425, F-91405 Orsay Cedex, France. E-mail address:
Sana.Louhichi@math.u-psud.fr}}

\date{December 7, 2009}

\maketitle

\begin{abstract}
\noindent In this article, we consider a %row-wise
stationary array $(X_{j,n})_{1 \leq j \leq n, n \geq 1}$ of random
variables with values in $\bR \verb2\2 \{0\}$ (which satisfy some
asymptotic dependence conditions), and the corresponding sequence
$(N_{n})_{n\geq 1}$ of point processes, where $N_{n}$ has the points
$X_{j,n}, 1\leq j \leq n$. Our main result identifies some explicit
conditions for the convergence of the sequence $(N_{n})_{n \geq 1}$,
in terms of the probabilistic behavior of the variables in the
array.
\end{abstract}


%60E07=infinitely divisible distributions
%60G55=point processes
%60G10=stationary processes
%60G57 Random measures


\section{Introduction}

The study of the asymptotic behavior of the sum (or the maximum) of
the row variables in an array $(X_{j,n})_{1 \leq j \leq n, n \geq
1}$ is one of oldest problem in probability theory. When the
variables are independent on each row, classical results identify
the limit to have an infinitely divisible distribution in the case
of the sum (see \cite{gnedenko-kolmogorov54}), and a max-infinitely
divisible distribution, in the case of the maximum (see
\cite{balkema-resnick77}). A crucial observation, which can be
traced back to \cite{resnick75}, \cite{weissman76} (in the case of
the maximum), and \cite{resnick-greenwood79} (in the case of the
sum) is that these results are deeply connected to the convergence
in distribution of the sequence $N_n=\sum_{j=1}^{n}\delta_{X_{j,n}},
n\geq 1$ of point processes to a Poisson process $N$. (See Section
5.3 of \cite{resnick87} and Section 7.2 of \cite{resnick07}, for a
modern account on this subject.)

Subsequent investigations showed that a similar connection exists in
the case of arrays which possess a row-wise dependence structure
(e.g. \cite{durrett-resnick78}). The most interesting case arises
when $X_{i,n}=X_i/a_n$, where $(X_i)_{i \geq 1}$ is a (dependent)
stationary sequence with regularly varying tails and $(a_n)_n$ is a
sequence of real numbers such that $n P(|X_1|>a_n) \to 1$ (see
\cite{DH95} and the references therein). In the dependent case, the
limit $N$ may not be a Poisson process, but belongs to the class of
infinitely divisible point processes (under generally weak
assumptions). These findings reveal that the separate study of the
point process convergence is an important topic, which may yield new
asymptotic results for triangular arrays.
%and give a new perspective to classical results.

In the present article, we consider an array $(X_{j,n})_{1 \leq j
\leq n,n \geq 1}$ whose row variables are asymptotically
independent, in the sense that the block $(X_{1,n}, \ldots,
X_{n,n})$ behaves asymptotically as $k_n$``smaller'' i.i.d. blocks,
a small block having the same distribution as $(X_{1,n}, \ldots,
X_{r_n,n})$, with $n \sim r_n k_n$. This condition, that we call
here (AD-1), was considered by many authors (e.g.
\cite{jakubowski93}
 \cite{jakubowski97}, \cite{DH95}, \cite{HHL88}, \cite{BL09}).

The rows of the array also possess an ``anti-clustering'' property
(AC), which specifies the dependence structure within a small block.
Intuitively, under (AC), it becomes improbable to find two points
$X_{j,n},X_{k,n}$ whose indices $j,k$ are situated in the same small
block at a distance larger than a fixed value $m$, and whose values
(in modulus) exceed a fixed threshold $\varepsilon>0$. Condition
(AC) appeared, in various forms, in the literature related to the
asymptotic behavior of the maximum  (e.g. \cite{leadbetter83},
\cite{obrien74}) or the sum (e.g. \cite{dedecker-louhichi05a},
\cite{dedecker-louhichi05b}, \cite{BJMW09}). In addition, we assume
the usual asymptotic negligibility (AN) condition for $X_{1,n}$.

Our main result says that under (AD-1), (AC) and (AN), the
convergence $N_n \stackrel{d}{\to} N$, where $N$ is an infinitely
divisible point process, reduces to the convergence of:
\begin{equation}
\label{cond-a-intro} n P(\max_{1 \leq j \leq m-1}|X_{j,n}| \leq x,
X_{m,n}>x), \ \mbox{and} \end{equation}
\begin{equation}
\label{cond-b-intro} n [P(A_{m,n}, \max_{1 \leq j \leq
m}|X_{j,n}|>x)- P(A_{m-1,n},\max_{1 \leq j \leq m-1}|X_{j,n}|>x)],
\end{equation}
 where
 $A_{m,n}$ is the event that at least $k_i$ among $X_{1,n}, \ldots,
 X_{m,n}$ lie in $B_i$, for all $i=1, \ldots,d$ (for
arbitrary $d,k_1, \ldots, k_d \in \bN$ and compact sets $B_1,
\ldots, B_d$).

The novelty of this result compared to the existing results (e.g.
Theorem 2.6 of \cite{BL09}), is the fact that the quantities
appearing in (\ref{cond-a-intro}) and (\ref{cond-b-intro}) speak
{\em explicitly} about the probabilistic behavior of the variables
in the array.

The article is organized as follows. In Section 2, we give the
statements of the main result (Theorem \ref{main-th}) and a
preliminary result (Theorem \ref{DH-th}). Section 3 is dedicated to
the proof of these two results. Section 4 contains a separate result
about the extremal index of a stationary sequence, whose proof is
related to some of the methods presented in this article.

\section{The Main Results}

We begin by introducing the terminology and the notation. Our main
reference is \cite{kall83}. We denote $\bR_{+}=[0,\infty)$,
$\bZ_{+}=\{0,1,2,\ldots\}$ and $\bN=\{1,2,\ldots\}$.

If $E$ is a locally compact Hausdorff space with a countable basis
(LCCB), we let $\cB$ be the class of all relatively compact Borel
sets in $E$, and $C_{K}^{+}(E)$ be the class of continuous functions
$f:E \to \bR_{+}$ with compact support. We let $M_p(E)$ be the class
of Radon measures on $E$ with values in $\bZ_+$ (endowed with the
topology of vague convergence), and $\cM_p(E)$ be the associated
Borel $\sigma$-field. For $\mu \in M_p(E)$ and $f \in C_K^{+}(E)$,
we denote $\mu(f)=\int_{E}f(x)\mu(dx)$. We denote by $o$ the null
measure.

Let $(\Omega, \cK,P)$ be a probability space. A measurable map $N:
\Omega \to M_p(E)$ is called a point process. Its distribution $P
\circ N^{-1}$ is determined by the Laplace functional
$L_N(f)=E(e^{-N(f)}), f \in C_K^{+}(E)$.

A point process $N$ is {\em infinitely divisible} if for any $k \geq
1$, there exist some i.i.d. point processes $(N_{i,k})_{1 \leq i
\leq k}$ such that $N \stackrel{d}{=} \sum_{i=1}^{k}N_{i,k}$. By
Theorem 6.1 of \cite{kall83}, the Laplace functional of an
infinitely divisible point process is given by:
$$L_{N}(f)=\exp \left\{-\int_{M_p(E) \verb2\2 \{o\}} (1-e^{-\mu(f)})
\lambda(d\mu) \right\}, \quad \forall f \in C_{K}^+(E),$$ where
$\lambda$ is a measure on $M_p(E) \verb2\2 \{o\}$, called the
canonical measure of $N$.

All the point processes considered in this article have their points
in $\bR \verb2\2 \{0\}$. For technical reasons, we embed $\bR
\verb2\2 \{0\}$ into the space $E=[-\infty,\infty] \verb2\2 \{0\}$.
Let $\cB$ be the class of relatively compact sets in $E$. Note that
$$[-x,x]^c:=[-\infty,-x) \cup (x,\infty] \in \cB, \quad \mbox{for
all} \ x>0.$$

We consider a triangular array $(X_{j,n})_{j \leq n, n \geq 1}$ of
random variables with values in $\bR \verb2\2 \{0\}$, such that
$(X_{j,n})_{j \leq n}$ is a strictly stationary sequence, for any $n
\geq 1$.

\begin{definition}
{\rm The triangular array $(X_{j,n})_{1 \leq j \leq n, n \geq 1}$
satisfies:

{\em (i)} {\bf condition (AN)} if
$$\limsup_{n \to \infty}nP(|X_{1,n}|
>\varepsilon )<\infty,\quad \mbox{for all} \ \varepsilon >0.$$

{\em (ii)} {\bf condition (AD-1)} if there exists $(r_n)_n \subset
\bN$ with $r_n \to \infty$ and $k_n=[n/r_n] \to \infty$, such that:
$$\lim_{n\to\infty}\left|E \left( e^{-\sum_{j=1}^{n} f(X_{j,n})}
\right)- \left\{E\left(e^{-\sum_{j=1}^{r_n}f(X_{j,n})}\right)
\right\}^{k_n} \right|=0, \quad \mbox{for all} \ f \in
C_{K}^{+}(E).$$

{\em (iii)} {\bf condition (AC)} if there exists $(r_n)_n \subset
\bN$ with $r_n \to \infty$, such that:
$$\lim_{m \to m_0}\limsup_{n \to \infty}n\sum_{j=m+1}^{r_n}
P(|X_{1,n}|>\varepsilon, |X_{j,n}|>\varepsilon)=0, \quad \mbox{for
all} \ \varepsilon >0,$$ where $m_0:=\inf \{m \in \bZ_+; \lim_{n \to
\infty}n\sum_{j=m+1}^{r_n} P(|X_{1,n}|>\varepsilon,
|X_{j,n}|>\varepsilon)=0, \ \mbox{for all} \ \varepsilon >0\}$. We
use the conventions: $\inf \emptyset=\infty$ and $\lim_{m \to
m_0}\phi(m)=\phi(m_0)$ if $m_0<\infty$.}
\end{definition}

\begin{remark}
{\rm (i) For each $n \geq 1$, let
$N_n=\sum_{j=1}^{n}\delta_{X_{j,n}}$ and $\tilde
N_n=\sum_{i=1}^{k_n} \tilde N_{i,n}$, where $(\tilde N_{i,n})_{i
\leq k_n}$ are i.i.d. copies of
$N_{r_n,n}=\sum_{j=1}^{r_n}\delta_{X_{j,n}}$. Under (AD-1),
$(N_n)_{n}$ converges in distribution if and only if $(\tilde
N_n)_n$ does, and the limits are the same.

(ii) Condition (AN) is an asymptotic negligibility condition which
ensures that $(\tilde N_{i,n})_{i \leq k_n, n \geq 1}$ is a {\em
null-array} of point processes, i.e. $P(\tilde N_{1,n}(B)>0) \to 0$
for all $B \in \cB$. By Theorem 6.1 of \cite{kall83} $\tilde N_n
\stackrel{d}{\to}  N$ if and only if
$$k_n E(1-e^{-N_{r_n,n}(f)} )
\to \int_{M_p(E) \verb2\2 \{o\}}(1-e^{-\mu(f)})\lambda(d\mu), \quad
\forall f \in C_{K}^{+}(E).$$ In this case, $N$ is an infinitely
divisible point process with canonical measure $\lambda$. }
\end{remark}

\begin{remark}
{\rm
 (i) Condition (AD-1) is satisfied by arrays whose row-wise dependence structure is of
mixing type (see e.g. Lemma 5.1 of \cite{BL09}).

(ii) Condition (AC) is satisfied with $m_0=m$ if $(X_{j,n})_{1 \leq
j \leq n}$ is $m$-dependent.

(iii) When $X_{j,n}=X_j/u_n$ and $m_0=1$, condition (AC) is known in
the literature as Leadbetter's condition $D'(\{u_n\})$ (see
\cite{leadbetter83}).

(iii) A condition similar to (AC) was used in \cite{BJMW09},
\cite{dedecker-louhichi05a} and \cite{dedecker-louhichi05b} for
obtaining the convergence of the partial sum sequence to an
infinitely divisible random variable (with finite variance).

}
\end{remark}

As in \cite{DH95}, let $M_0 = \{\mu \in M_p(\bR \verb2\2 \{0\}); \mu
\not = o, \ \exists \ x \in (0,\infty) \ \mbox{such that} \
\mu([-x,x]^c)=0\}$. If $\mu=\sum_{j \geq 1} \delta_{t_j} \in M_0$,
we let $x_{\mu}:=\sup_{j \geq 1}|t_j|<\infty$. For each $x>0$, let
$$M_x =\{\mu \in M_0; \mu([-x,x]^c)>0\}=\{\mu \in M_0; x_{\mu}>x\}.$$

%Note that $M_0$ and $M_x,x>0$ are closed sets of $M_p(E)$ with
%respect to the vague topology (see Lemma 3.3 of \cite{BL09p}).

%We define the measure $\nu$ on $(0,\infty)$ by:
%$$\nu^*(x,\infty)=\lambda(M_x) \quad \mbox{for all} \ x>0.$$

Recall that $x$ is a {\em fixed atom} of a point process $N$ if
$P(N\{x\}>0)>0$. To simplify the writing, we introduce some
additional notation. If $x>0$ and $\lambda$ is a measure on $M_p(E)$
with $\lambda(M_0^c)=0$, we let
$$\cB_{x,\lambda}=\{B \in \cB;\lambda(\{\mu \in M_x; \mu(\partial
B)>0\})=0\},$$ and $\cJ_{x,\lambda}$ be the class of sets
$M=\cap_{i=1}^{d}\{\mu \in M_p(E); \mu(B_i) \geq k_i\}$ for some
$B_i \in \cB_{x,\lambda}$, $k_i \geq 1$ (integers) and $d \geq 1$.

The following result is a refinement of Theorem 3.6 of \cite{BL09p}.

\begin{theorem}
\label{DH-th} Suppose that $(X_{j,n})_{1 \leq j \leq n, n \geq 1}$
satisfies (AN) and (AD-1) (with sequences $(r_n)_n$ and $(k_n)_n$).
Let $N$ be an infinitely divisible point process on $\bR \verb2\2
\{0\}$ with canonical measure $\lambda$. Let $D$ be the set of fixed
atoms of $N$ and $D'=\{x>0; x \in D \ \mbox{or} \ -x \in D\}$.

The following statements are equivalent:

(i) $N_n \stackrel{d}{\to} N$;

(ii) We have $\lambda(M_0^c)=0$, and the following two conditions
hold:
\begin{eqnarray*}
(a) & & k_n P(\max_{j \leq r_n}|X_{j,n}|>x)  \to \lambda(M_x), \
\mbox{for any} \
x>0, x \not \in D', \\
(b) & & k_n P(N_{r_n,n} \in M, \max_{j \leq r_n}|X_{j,n}|>x) \to
\lambda(M \cap M_x), \ \mbox{for any} \ x>0, x \not \in
D' \\
& & \mbox{and for any set} \ M \in \cJ_{x,\lambda}.
\end{eqnarray*}
\end{theorem}

%\begin{lemma}
%If $N_n \stackrel{d}{\to} N$, where $N$ is an infinitely divisible
%point process of canonical measure $\lambda$, then
%$\lambda(M_0^c)=0$ and $M_n =\max_{1 \leq n}|X_{j,n}|
%\stackrel{d}{\to} Y$, where $P(Y \leq x)=e^{-\nu(x,\infty)}$ for any
%$x>0$.
%\end{lemma}

%\noindent {\bf Proof:} To be written.

%\vspace{30mm}

For each $1 \leq m \leq n$, let
$N_{m,n}=\sum_{j=1}^{m}\delta_{X_{j,n}}$ and $M_{m,n}=\max_{j \leq
m}|X_{j,n}|$, with the convention that $M_{0,n}=0$. The next theorem is the main result of this article,
and gives an explicit form for conditions (a) and (b), under the
additional anti-clustering condition (AC).

\begin{theorem}
\label{main-th} Let $(X_{j,n})_{1 \leq j \leq n, n \geq 1}$ and $N$
be as in Theorem \ref{DH-th}. Suppose in addition that (AC) holds,
with the same sequence $(r_n)_n$ as in (AD-1).

The following statements are equivalent:

(i) $N_n \stackrel{d}{\to} N$;

(ii) We have $\lambda(M_0^c)=0$ and the following two conditions
hold:
\begin{eqnarray*}
(a') & & \lim_{m \to m_0} \limsup_{n \to \infty}
|n[P(M_{m,n}>x)-P(M_{m-1,n}>x)]-\lambda(M_x)| =0, \ \mbox{for any}
\\ & &
 x>0, x \not \in D', \\
(b') & & \lim_{m \to m_0} \limsup_{n \to \infty} |n[P(N_{m,n} \in M,
M_{m,n}>x)-P(N_{m-1,n} \in M, M_{m-1,n}>x)]  \\
& & -\lambda(M \cap M_x)|=0, \ \mbox{for any} \ x>0, x \not \in D' \
\mbox{and for any set} \ M \in \cJ_{x,\lambda}.
\end{eqnarray*}
\end{theorem}

\begin{remark}
{\rm Note that $$P(M_{m,n}>x)-P(M_{m-1,n}>x)=P(\max_{1 \leq j \leq
m-1}|X_{j,n}| \leq x,|X_{m,n}|>x).$$
 }
\end{remark}

\begin{remark}
{\rm Suppose that $m_0=1$ in Theorem \ref{main-th}. One can prove that in this case, the limit $N$ is a Poisson process of intensity $\nu$ given by:
$$\nu(B)=\lambda(\{\mu\in M_p(E) \verb2\2 \{o\}; \mu(B)=1\}), \quad \forall B \in \cB.$$
}
\end{remark}

\section{The Proofs}

\subsection{Proof of Theorem \ref{DH-th}}

Before giving the proof, we need some preliminary results.

\begin{lemma}
\label{boundary-M} Let $E$ be a LCCB space and
$M=\cap_{i=1}^{d}\{\mu \in M_p(E); \mu(B_i) \geq k_i\}$ for
some $B_i \in \cB$, $k_i \geq 1$ (integers) and $d \geq 1$. Then:\\
(i)  $M$ is closed (with respect to the vague topology); \\
(ii) $\partial M \subset \cup_{i=1}^{d}\{\mu \in M_p(E);
\mu(\partial B_i)>0\}$.
\end{lemma}

\noindent {\bf Proof:} Note that $\partial M \subset
\cup_{i=1}^{d}\partial M_i$, where $M_i=\{\mu \in M_p(E); \mu(B_i)
\geq k_i \}$. Since the finite intersection of closed sets is a
closed set, it suffices to consider the case $d=1$, i.e. $M=\{\mu
\in M_p(E); \mu(B) \geq k\}$ for some $B \in \cB$ and $k \geq 1$.

(i) Let $(\mu_n)_n \subset M$ be such that $\mu_n
\stackrel{v}{\to}\mu$. If $\mu(\partial B)=0$, then $\mu_n(B) \to
\mu(B)$, and since $\mu_n(B) \geq k$ for all $n$, it follows that
$\mu(B) \geq k$. If not, we proceed as in the proof of Lemma 3.15 of
\cite{resnick87}. Let $B^{\delta}$ be a $\delta$-swelling of $B$.
Then $S=\{\delta \in (0,\delta_0]; \mu(\partial B^{\delta})>0 \}$ is
a countable set. By the previous argument, $\mu(B^{\delta}) \geq k$
for all $\delta \in (0,\delta_0] \verb2\2 S$. Let $(\delta_n)_n \in
(0,\delta_0] \verb2\2 S$ be such that $\delta_n \downarrow 0$. Since
$\mu(B^{\delta_n}) \geq k$ for all $n$, and $\mu(B^{\delta_n})
\downarrow \mu(B)$, it follows that $\mu(B) \geq k$, i.e. $\mu \in
M$.

(ii) By part (i), $\partial M=\bar M \verb2\2 M^o= M \cap (M^o)^c$.
We will prove that $\partial M \subset \{\mu \in M; \mu(\partial
B)>0\}$, or equivalently
$$A:=\{\mu \in M; \mu(\partial B)=0\} \subset M^o.$$

Since $M^o$ is the largest open set included in $M$ and $A \subset
M$, it suffices to show that $A$ is open. Let $\mu \in A$ and
$(\mu_n)_n \subset M_p(E)$ be such that $\mu_n
\stackrel{v}{\to}\mu$. Then $\mu_{n}(B) \to \mu(B)$, and since
$\mu(B) \geq k$ and $\{\mu_n(B)\}_n$ are integers, it follows that
$\mu_n(B) \geq k$ for all $n \geq N_1$, for some $N_1$.

On the other hand, $\mu_n(\partial B) \to \mu(\partial B)$, since
$\partial B \in \cB$ and $\mu(\partial B)=0$ (note that $\partial
(\partial B)=\partial B$). Since $\mu(\partial B)=0$ and
$\{\mu_n(\partial B)\}_n$ are integers, it follows that
$\mu_n(\partial B)=0$ for all $n \geq n_2$, for some $n_2$. Hence
$\mu_n \in A$ for all $n \geq \max\{n_1,n_2\}$. $\Box$

\begin{lemma}
\label{weak-conv-proba}Let $E$ be a LCCB space and $(Q_n)_n,Q$ be
probability measures on $M_p(E)$. Let $\cB_{Q}$ be the class of all
sets $B \in \cB$ which satisfy: $$Q(\{\mu \in M_p(E) ; \mu(\partial
B)>0\})=0,$$ and $\cJ_{Q}$ be the class of sets
$M=\cap_{i=1}^{d}\{\mu \in M_p(E); \mu(B_i) \geq k_i \}$ for some
$B_i \in \cB_{Q}$, $k_i \geq 1$ (integers) and $d \geq 1$.

Then $Q_n \stackrel{w}{\to} Q$ if and only if $Q_n(M) \to Q(M)$ for
all $M \in \cJ_{Q}$.
\end{lemma}

\noindent {\bf Proof:} Let $(N_n)_n, N$ be point processes on $E$,
defined on a probability space $(\Omega,\cF,P)$, such that $P \circ
N_n^{-1}=Q_n$ for all $n$, and $P \circ N^{-1}=Q$. Note that
$\cB_{Q}=\cB_{N}:=\{B \in \cB; N(\partial B)=0 \ \mbox{a.s.} \}$.

By definition, $N_n \stackrel{d}{\to} N$ if and only if $Q_n
\stackrel{w}{\to} Q$. By Theorem 4.2 of \cite{kall83}, $N_n
\stackrel{d}{\to} N$ if and only if
$$(N_n(B_1), \ldots, N_n(B_d)) \stackrel{d}{\to} (N(B_1), \ldots,
N(B_d))$$ for any $B_1, \ldots, B_d \in \cB_{N}$ and for any $d \geq
1$. Since these random vectors have values in $\bZ_{+}^{d}$, the
previous convergence in distribution is equivalent to:
$$P(N_n(B_1) =k_1, \ldots, N_n(B_d) = k_d) \to P(N(B_1) = k_1, \ldots,
N(B_d) = k_d)$$ for any $k_1, \ldots, k_d \in \bZ_{+}$, which is in
turn equivalent to
$$P(N_n(B_1) \geq k_1, \ldots, N_n(B_d) \geq k_d) \to P(N(B_1) \geq k_1, \ldots,
N(B_d) \geq k_d)$$ for any $k_1, \ldots, k_d \in \bZ_{+}$. Finally,
it suffices to consider only integers $k_i \geq 1$ since, if there
exists a set $I \subset \{1,\ldots,d\}$ such that $k_i=0$ for all $i
\in I$ and $k_i \geq 1$ for $i \not \in I$, then $P(N_n(B_1) \geq
k_1, \ldots, N_n(B_d) \geq k_d)=P(N_n(B_i) \geq k_i, i \not \in I)
\to P(N(B_i) \geq k_i,i \not \in I )=P(N(B_1) \geq k_1, \ldots,
N(B_d) \geq k_d)$. $\Box$

\vspace{3mm}

\noindent {\bf Proof of Theorem \ref{DH-th}:} Note that $\{\max_{j
\leq r_n}|X_{j,n}|>x\}=\{N_{r_n,n} \in M_x\}$.

{\em Suppose that (i) holds.} As in the proof of Theorem 3.6 of
\cite{BL09p}, it follows that $\lambda(M_0^c)=0$ and (a) holds.
Moreover, we have $P_{n,x} \stackrel{w}{\to}P_{x}$ where $P_{n,x}$
and $P_{x}$ are probability measures on $M_p(E)$ defined by:
$$P_{n,x}(M)=\frac{k_nP(N_{r_n,n} \in M \cap M_x)}{k_n P(N_{r_n,n} \in
M_x)} \quad \mbox{and} \quad P_x(M)=\frac{\lambda(M \cap
M_x)}{\lambda(M_x)}.$$

\noindent Therefore, $P_{n,x}(M) \to P_x(M)$ for any $M \in
\cM_p(E)$ with $P_{x}(\partial M)=0$. Since $k_n P(N_{r_n,n} \in
M_x) \to \lambda(M_x)$ (by (a)), it follows that
\begin{equation}
\label{conv-b} k_nP(N_{r_n,n} \in M \cap M_x) \to \lambda(M \cap
M_x),
\end{equation}
for any $M \in \cM_p(E)$ with $\lambda(\partial M \cap M_x)=0$.

In particular, (\ref{conv-b}) holds for a set $M=\cap_{i=1}^{d}\{\mu
\in M_p(E); \mu(B_i) \geq 1\}$, with $B_i \in \cB_{x,\lambda}$, $k_i
\geq 1$ (integers) and $d \geq 1$. To see this, note that by Lemma
\ref{boundary-M}, $\partial M \cap M_x \subset \cup_{i=1}^{d} \{\mu
\in M_x; \mu(\partial B_i)>0\}$, and hence
$$\lambda(\partial M \cap M_x) \leq \sum_{i=1}^{d}\lambda (\{ \mu \in M_x;
\mu(\partial B_i)>0\})=0.$$

{\em Suppose that (ii) holds.} As in the proof of Theorem 3.6 of
\cite{BL09p}, it suffices to show that $P_{n,x} \stackrel{w}{\to}
P_{x}$. This follows by Lemma \ref{weak-conv-proba}, since the class
of sets $B \in \cB$ which satisfy:
$$P_x(\{\mu \in M_p(E); \mu(\partial B)>0 \})=0$$
coincides with $\cB_{x,\lambda}$. $\Box$

\subsection{Proof of Theorem \ref{main-th}}

We begin with an auxiliary result, which is of independent interest.

\begin{lemma}
\label{m-block-lemma} Let $h:\bR^d \to \bR$ be a twice continuously
differentiable function, such that
\begin{equation}
\label{cond-der-h} \|D^2 h \|_{\infty}:=\max_{i,j =1,\ldots,d}
\sup_{{\bf x} \in \bR^d}\left|\frac{\partial^2 h}{\partial x_i
\partial x_j}({\bf x})\right|<\infty.
\end{equation}
Let $({\bf Y}_i)_{i \geq 1}$ be a strictly stationary sequence of
$d$-dimensional random vectors with ${\bf Y}_i=(Y_{i}^{(1)}, \ldots,
Y_{i}^{(d)})$. Let ${\bf S}_n=\sum_{i=1}^{n}{\bf Y}_i$ for $n \geq
1$ and ${\bf S}_0={\bf 0}$. Then for any $1 \leq m \leq r$,
$$|E[h({\bf S}_r)]-r E[h({\bf S}_m)-h({\bf
S}_{m-1})]| \leq m |E[h({\bf S}_m)]+E[h({\bf S}_{m-1})]|+$$
$$\|D^2 h \|_{\infty} \sum_{k=0}^{r_n-m}
\sum_{i,l=1}^{d}E|S_k^{(i)}Y_{k+m}^{(l)}|.$$
\end{lemma}

\noindent {\bf Proof:} As in Lemma 3.2 of \cite{jakubowski97} (see
also Theorem 2.6 of \cite{BL09}), we have:
\begin{eqnarray*}
E[h({\bf S}_r)] &=& E[h({\bf S}_{m-1})]+\sum_{k=0}^{r-m}E[h({\bf
S}_{k+m})-h({\bf S}_{k+m-1})] \\
r E[h({\bf S}_m)-h({\bf S}_{m-1})] &=& (m-1) E[h({\bf S}_m)-h({\bf
S}_{m-1})] + \\
& & \sum_{k=0}^{r-m} E[h({\bf S}_{k+m}-{\bf S}_k) -h({\bf
S}_{k+m-1}-{\bf S}_k)],
\end{eqnarray*}
where the second equality is due to the stationarity of $({\bf
Y}_i)_{i}$. Taking the difference, we get:
\begin{eqnarray*}
\lefteqn{E[h({\bf S}_r)]- r E[h({\bf S}_m)-h({\bf S}_{m-1})] = m
E[h({\bf S}_{m-1})] -(m-1) E[h({\bf S}_m)] + } \\
& & \sum_{k=0}^{r-m}E \{[h({\bf S}_{k+m})-h({\bf S}_{k+m}-{\bf
S}_{k})]-[h({\bf S}_{k+m-1})-h({\bf S}_{k+m-1}-{\bf S}_{k})]
\}=:I_1+I_2.
\end{eqnarray*}

Clearly $|I_1| \leq  m |E[h({\bf S}_{m-1})] + E[h({\bf S}_m)]|$. For
treating $I_2$, we use the Taylor's formula (with integral
remainder) for twice continuously differentiable functions $f:\bR^d
\to \bR$:
\begin{equation}
\label{taylor} f(x)-f(x_0)=\sum_{i=1}^{d}(x^{(i)}-x_0^{(i)})
\int_0^1 \frac{\partial f}{\partial x_i}(x-s(x-x_0))ds.
\end{equation}
We get:
\begin{eqnarray*}
h({\bf S}_{k+m})-h({\bf S}_{k+m}-{\bf S}_k)&=&
\sum_{i=1}^{d}S_k^{(i)} \int_0^1 \frac{\partial h}{\partial
x_i}({\bf S}_{k+m}-x {\bf S}_{k})dx,  \\
h({\bf S}_{k+m-1})-h({\bf S}_{k+m-1}-{\bf S}_k)&=&
\sum_{i=1}^{d}S_k^{(i)} \int_0^1 \frac{\partial h}{\partial
x_i}({\bf S}_{k+m-1}-x {\bf S}_{k})dx.
\end{eqnarray*}

Taking the difference of the last two equations, and using
(\ref{taylor}) for $f=\partial h/\partial x_i$ with $i=1,\ldots,d$,
we obtain:
\begin{eqnarray*}
\lefteqn{[h({\bf S}_{k+m})-h({\bf S}_{k+m}-{\bf S}_k)]-[h({\bf
S}_{k+m-1})-h({\bf S}_{k+m-1}-{\bf S}_k)] = }\\
& & \sum_{i=1}^{d}S_k^{(i)} \int_0^1 \left[ \frac{\partial
h}{\partial x_i}({\bf S}_{k+m}-x {\bf S}_{k})-\frac{\partial
h}{\partial x_i}({\bf S}_{k+m-1}-x {\bf S}_{k}) \right]dx \\
&=& \sum_{i=1}^{d}S_k^{(i)} \int_0^1 \sum_{l=1}^{d} Y_{k+m}^{(l)}
\int_0^1 \frac{\partial^2 h}{\partial x_i \partial x_l}(({\bf
S}_{k+m}-x {\bf S}_k)-\theta {\bf Y}_{k+m})d\theta dx.
\end{eqnarray*}

\noindent From here we conclude that:
$$|[h({\bf S}_{k+m})-h({\bf S}_{k+m}-{\bf S}_k)]-[h({\bf
S}_{k+m-1})-h({\bf S}_{k+m-1}-{\bf S}_k)]| \leq \|D^2
h\|_{\infty}\sum_{i,l=1}^{d}|S_{k}^{(i)}Y_{k+m}^{(l)}|,$$ which
yields the desired estimate for $I_2$. $\Box$

\begin{proposition}
\label{difference} Let $E$ be a LCCB space. For each $n \geq 1$, let
$(X_{j,n})_{j \leq n}$ be a strictly stationary sequence of
$E$-valued random variables, such that:
\begin{equation}
\label{AN-arb-space} \limsup_{n \to \infty} n P(X_{1,n} \in
B)<\infty, \ \mbox{for all} \ B \in \cB.
\end{equation}

Suppose that there exists $(r_n)_n \subset \bN$ with $r_n \to
\infty$ and $k_n=[n/r_n] \to \infty$, such that:
\begin{equation}
\label{AC-arb-space} \lim_{m \to m_0} \limsup_{n \to \infty} n
\sum_{j=m+1}^{r_n}P(X_{1,n} \in B, X_{j,n} \in B)=0, \ \mbox{for
all} \ B \in \cB,
\end{equation}
where $m_0=:\{m \in \bZ_{+}; \lim_{n \to \infty} n
\sum_{j=m+1}^{r_n}P(X_{1,n} \in B, X_{j,n} \in B)=0,  \ \mbox{for
all} \ \ B \in \cB\}$. Let $N_{m,n}=\sum_{j=1}^{m}\delta_{X_{j,n}}$.
Then
$$\lim_{m \to m_0} \limsup_{n \to \infty} |k_n P(N_{r_n,n} \in M)-n
[P(N_{m,n} \in M)-P(N_{m-1,n} \in M)]|=0,$$ for any set
$M=\cap_{i=1}^{d} \{\mu \in M_p(E); \mu(B_i) \geq k_i \}$, with $B_i
\in \cB$, $k_i \geq 1$ (integers) and $d \geq 1$.
\end{proposition}

\noindent {\bf Proof:} Let $h:\bR_+^d \to \bR_{+}$ be a twice
continuously differentiable function which satisfies
(\ref{cond-der-h}), such that $h(x_{1}, \ldots,x_d) \leq x_1+ \ldots
+x_d$ for all $(x_1, \ldots,x_d) \in \bR_+^d$, and
$$h(x_{1},\ldots,x_d)=\left\{
\begin{array}{ll} 0 & \mbox{if $x_i \leq k_i-1$ for some $i=1,\ldots,d$} \\
1 & \mbox{if $x_i \geq k_i$ for all $i=1,\ldots,d$}
\end{array}  \right.$$
Note that:
\begin{equation}
\label{property-h} h(x_1,\ldots,x_d)=1_{\{x_1 \geq k_1,\ldots,x_d
\geq k_d\}} \quad \mbox{for all} \ x_1, \ldots, x_d \in \bZ_{+}.
\end{equation}

For any $n \geq 1$, we consider strictly stationary sequence of
$d$-dimensional random vectors $\{{\bf Y}_{j,n}=(Y_{j,n}^{(1)},
\ldots, Y_{j,n}^{(d)}), 1 \leq j \leq n\}$ defined by:
$$Y_{j,n}^{(i)}=1_{\{X_{j,n} \in B_i\}}, \quad \mbox{for any} \
i=1,\ldots,d.$$ Using (\ref{property-h}), we obtain for any $1 \leq
m \leq n$,
\begin{eqnarray*}
\lefteqn{P(N_{m,n} \in M)=P(N_{m,n}(B_1) \geq k_1, \ldots,
N_{m,n}(B_d) \geq k_d) =}\\
& & P(\sum_{j=1}^{m}Y_{j,n}^{(1)} \geq k_1, \ldots,
\sum_{j=1}^{m}Y_{j,n}^{(d)} \geq k_d ) = E[1_{\{\sum_{j=1}^{m}
Y_{j,n}^{(1)} \geq k_1, \ldots,
\sum_{j=1}^{m} Y_{j,n}^{(d)} \geq k_d \}}] =\\
 & & E[h(\sum_{j=1}^{m}Y_{j,n}^{(1)},\ldots,\sum_{j=1}^{m}Y_{j,n}^{(d)}
 )]=E[h(\sum_{j=1}^{m}{\bf Y}_{j,n})].
\end{eqnarray*}

Using Lemma \ref{m-block-lemma}, and letting $C=\|D^2 h\|_{\infty}$,
we obtain:
\begin{eqnarray}
\nonumber \lefteqn{ k_n |P(N_{r_n,n} \in M)-r_n[P(N_{m,n} \in
M)-P(N_{m-1,n} \in
M)]| \leq } \\
\nonumber & & m k_n \{E[h(\sum_{j=1}^{m}{\bf
Y}_{j,n})]+E[h(\sum_{j=1}^{m-1}{\bf Y}_{j,n})] \} + C k_n
\sum_{i,l=1}^{d}\sum_{k=0}^{r_n-m}E (\sum_{j=1}^k
Y_{j,n}^{(i)}Y_{k+m,n}^{(l)}) \\
\label{sum-I1-I2} & &=:I_{m,n}^{(1)}+C I_{m,n}^{(2)}
\end{eqnarray}

\noindent Using the fact that $h({\bf x}) \leq \sum_{i=1}^{d}x_i$,
and the stationary of $(X_{j,n})_{j \leq n}$,
\begin{eqnarray*}
I_{m,n}^{(1)} & \leq & 2m k_nE(\sum_{j=1}^{m}
\sum_{i=1}^{d}Y_{j,n}^{(i)})=2m k_n \sum_{i=1}^{d}
\sum_{j=1}^{m}P(X_{j,n} \in B_i) \\
&=&2m^2 k_n \sum_{i=1}^{d} P(X_{1,n} \in B_i)\leq 2m^2 \frac{1}{r_n}
\sum_{i=1}^{d} n P(X_{1,n} \in B_i).
\end{eqnarray*}

\noindent From (\ref{AN-arb-space}), it follows that $\lim_{n \to
\infty}I_{m,n}^{(1)}=0$ for all $m$, and hence
\begin{equation}
\label{conv-I1} \lim_{m \to m_0}\limsup_{n \to
\infty}I_{m,n}^{(1)}=0.
\end{equation}

Using the stationarity of $(X_{j,n})_{j \leq n}$, and letting
$B=\cup_{i=1}^{d}B_i \in \cB$,
\begin{eqnarray*}
I_{m,n}^{(2)} &=& k_n
\sum_{i,l=1}^{d}\sum_{k=0}^{r_n-m}\sum_{j=1}^{k}
P(X_{j,n} \in B_i, X_{k+m,n} \in  B_l) \\
&=& k_n \sum_{i,l=1}^{d} \sum_{j=m+1}^{r_n} (r_n-j+1) P(X_{1,n} \in
B_i, X_{j,n} \in B_l) \\
& \leq &  n \sum_{i,l=1}^{d}  \sum_{j=m+1}^{r_n} P(X_{1,n} \in B_i,
X_{j,n} \in B_l) \\
& \leq & d^2 n  \sum_{j=m+1}^{r_n} P(X_{1,n} \in B, X_{j,n} \in B).
\end{eqnarray*}

\noindent From (\ref{AC-arb-space}), it follows that:
\begin{equation}
\label{conv-I2} \lim_{m \to m_0}\limsup_{n \to
\infty}I_{m,n}^{(2)}=0.
\end{equation}

From, (\ref{sum-I1-I2}), (\ref{conv-I1}) and (\ref{conv-I2}), it
follows that:
$$\lim_{m \to m_0} \limsup_{n \to \infty} k_n |P(N_{r_n,n} \in M)-r_n[P(N_{m,n} \in
M)-P(N_{m-1,n} \in M)]|=0.$$

\noindent Note that $\lim_{n \to \infty} (n-k_n r_n)|P(N_{m,n} \in
M)-P(N_{m-1,n} \in M)|=0$ for all $m$, and hence
$$\lim_{m \to m_0} \limsup_{n \to \infty} (n-k_n r_n)|P(N_{m,n} \in
M)-P(N_{m-1,n} \in M)|=0.$$ The conclusion follows. $\Box$

\begin{corollary}
\label{cor-difference} For each $n \geq 1$, let $(X_{j,n})_{1 \leq j
\leq n}$ be a strictly stationary sequence of random variables with
values in $\bR \verb2\2 \{0\}$. Suppose that $(X_{j,n})_{1 \leq j
\leq n, n \geq 1}$ satisfies (AN) and (AC).

For any $1 \leq m \leq n$, let
$N_{m,n}=\sum_{j=1}^{m}\delta_{X_{j,n}}$ and $M_{m,n}=\max_{j \leq
m}|X_{j,n}|$. Then,
\begin{eqnarray*}
& & \lim_{m \to m_0} \limsup_{n \to \infty} |k_n P(M_{r_n,n}>x)-n
[P(M_{m,n}>x)-P(M_{m-1,n}>x)]|=0 \\
& & \lim_{m \to m_0} \limsup_{n \to \infty} |k_n P(N_{r_n,n} \in M,
M_{r_n,n}>x)-n [P(N_{m,n} \in M,M_{m,n}>x)-\\
& & P(N_{m-1,n} \in M,M_{m-1,n}>x)]|=0,
\end{eqnarray*}
for any $x>0$, and for any set $M=\cap_{i=1}^{d} \{\mu \in M_p(E);
\mu(B_i) \geq k_i \}$, with $B_i \in \cB$, $k_i \geq 1$ (integers)
and $d \geq 1$.
\end{corollary}

\noindent {\bf Proof:} Since $\{M_{m,n}>x\}=\{N_{m,n}([-x,x]^c) \geq
1 \}$ for any $1 \leq m \leq n$, the result follows from Proposition
\ref{difference}. $\Box$

\vspace{3mm}

\noindent {\bf Proof of Theorem \ref{main-th}:} The result follows
from Theorem \ref{DH-th} and Corollary \ref{cor-difference}. $\Box$

\section{The extremal index}

In this section, we give a recipe for calculating the extremal index
of a stationary sequence, using a method which is similar to that
used for proving Theorem \ref{main-th}, in a simplified context.
Although this recipe (given by Theorem \ref{extremal-th} below)
seems to be known in the literature (see \cite{obrien74},
\cite{obrien87}, \cite{smith92}), we decided to include it here,
since we could not find a direct reference for its proof.

We recall the following definition.

\begin{definition}
{\rm Let $(X_{j})_{j \geq 1}$ be a strictly stationary sequence of
random variables. {\bf The extremal index} of the sequence $(X_j)_{j
\geq 1}$, if it exists, is a real number $\theta$ with the following
property: for any $\tau>0$, there exists a sequence
$(u_n^{(\tau)})_n\subset \bR$ such that $n P(X_1>u_{n}^{(\tau)}) \to
\tau$ and $P(\max_{j \leq n}X_j \leq u_n^{(\tau)}) \to e^{-\tau
\theta}$. }
\end{definition}

\noindent In particular, for $\tau=1$, we denote $u_n^{(1)}=u_n$,
and we have
\begin{equation}
\label{def-u-n} n P(X_1>u_{n})  \to 1 \  \quad \mbox{and} \quad
P(\max_{j \leq n}X_j \leq u_n) \to e^{- \theta}.
\end{equation}
 It is clear that if it exists, $\theta
\in [0,1]$.

\begin{remark}
{\rm The extremal index of an i.i.d sequence exists and is equal to
$1$. }
%(ii)  The extremal index of a negatively associated sequence, {\em
%if it exists}, is equal to $1$. (Recall that a sequence $(X_j)_j$ is
%negatively associated if $${\rm Cov}(f(X_i,i \in I), g(X_i,i \in J))
%\leq 0$$ for any disjoint sets $I,J$, and for any coordinate-wise
%non-decreasing functions $f,g$.) }
\end{remark}

The following definition was originally considered in
\cite{obrien74}.

\begin{definition}
{\rm We say that $(X_j)_{j \geq 1}$ satisfies {\bf condition (AIM)}
(or admits an asymptotic independence representation for the
maximum) if there exists $(r_n)_n \subset \bN$ with $r_n \to \infty$
and $k_n=[n/r_n]\to \infty$, such that:
$$\left|P(\max_{j \leq n} X_j \leq u_n) -P(\max_{j \leq r_n}X_j \leq u_n)^{k_n} \right|
\to 0.$$ }
\end{definition}

\begin{remark}
{\rm It is known that (Leadbetter's) condition $D(\{u_n\})$ implies
(AIM) (see Lemma 2.1 of \cite{leadbetter83}). Recall that
$(\xi_j)_j$ satisfies condition $D(\{u_n\})$ if there exists a
sequence $(m_n)_n \subset \bN$, such that $m_n=o(n)$ and
$\alpha_n(m_n)\to 0$, where
$$\alpha_n(m)=\sup_{I,J} |P(\max_{j \in I}X_j \leq u_n,\max_{j \in J}X_j \leq u_n)
-P(\max_{j \in I}X_j \leq u_n)P(\max_{j \in J}X_j \leq u_n)|,$$
where the supremum ranges over all disjoint subsets $I,J$ of $\{1,
\ldots,n\}$, which are separated by a block of length greater of
equal than $m$.

}
\end{remark}

The following theorem is the main result of this section.

\begin{theorem}
\label{extremal-th} Let $(X_{j})_{j \geq 1}$ be a strictly
stationary sequence whose extremal index $\theta$ exists, and
$(u_n)_n$ be a sequence of real numbers satisfying (\ref{def-u-n}).

Suppose that $(X_j)_{j \geq 1}$ satisfies (AIM), and in addition,
\begin{equation}
\label{AC-extremal} \lim_{m \to m_0} \limsup_{n \to \infty}n
\sum_{j=m+1}^{r_n}P(X_1>u_n, X_{j}>u_n)=0,
\end{equation}
where $m_0:=\inf \{m \in \bZ_{+}; \lim_{n \to \infty}n
\sum_{j=m+1}^{r_n}P(X_1>u_n, X_{j}>u_n)=0\}$.

Then
\begin{equation}
\label{new-def-theta} \theta=\lim_{m \to m_0} \limsup_{n \to \infty}
n P(\max_{1 \leq j \leq m-1}X_{j} \leq u_n,X_{m}>u_n).
\end{equation}
\end{theorem}

Due to the stationarity, and the fact that $nP(X_{1}>u_n) \to 1$,
(\ref{new-def-theta}) can be written as:
\begin{eqnarray*}
\theta &=& \lim_{m \to m_0} \limsup_{n \to \infty} n P(\max_{2 \leq
j \leq m}X_{j} \leq u_n,X_{1}>u_n) \\
&=& \lim_{m \to m_0} \limsup_{n \to \infty} P(\max_{2 \leq j \leq
m}X_{j} \leq u_n|X_{1}>u_n),
\end{eqnarray*}
which coincides with (2.3) of \cite{smith92}.
\begin{remark}{\rm Let $(Y_i)_{i\geq 1}$ be a sequence of i.i.d. random variables and $X_i=\max(Y_{i},\cdots, Y_{i+m-1}).$
Then $(X_i)_{i\geq 1}$ satisfies condition (\ref{AC-extremal}), since it is an $m$-dependent sequence. A direct calculation shows that the extremal index of $(X_i)_{i\geq 1}$ exists
and is equal to $1/m$, which can be deduced also from (\ref{new-def-theta}).}
\end{remark}

The proof of Theorem \ref{extremal-th} is based on some intermediate
results.

\begin{proposition}
\label{calcul-theta} Let $(X_{j})_{j \geq 1}$ be a strictly
stationary sequence whose extremal index $\theta$ exists, and
$(u_n)_n$ be a sequence of real numbers satisfying (\ref{def-u-n}).
If $(X_j)_{j \geq 1}$ satisfies (AIM), then $k_n P(\max_{j \leq
r_n}X_j>u_n) \to \theta$.
\end{proposition}

\noindent {\bf Proof:} Due to (AIM), $P(\max_{j \leq r_n}X_j \leq
u_n)^{k_n} \to e^{-\theta}$. The result follows, since
$$P(\max_{j \leq r_n}X_j \leq u_n)^{k_n}=
\left(1-\frac{k_nP(\max_{j \leq r_n}X_j>u_n)}{k_n} \right)^{k_n}.$$
$\Box$

\begin{proposition}
\label{difference2} Let $(X_j)_{j \geq 1}$ be a strictly stationary
sequence such that:
$$\limsup_{n \to \infty}n P(X_1>u_n)<\infty.$$
Suppose that there exists $(r_n)_n \subset \bN$ with $r_n \to
\infty$ and $k_n=[n/r_n] \to \infty$, such that (\ref{AC-extremal})
holds. Then
$$\lim_{m \to m_0} \limsup_{n \to \infty} |k_n P(\max_{j \leq r_n}X_j>u_n)-n
[P(\max_{j \leq m}X_j>u_n)-P(\max_{j \leq m-1}X_j>u_n)]|=0.$$
\end{proposition}

\noindent {\bf Proof:} The argument is the same as in the proof of
Proposition \ref{difference}, using Lemma \ref{m-block-lemma}. More
precisely, we let $h:\bR_{+} \to \bR_{+}$ be a twice continuously
differentiable such that $\|h''\|_{\infty}<\infty$, $h(0)=0$,
$h(1)=1$ if $y \geq 1$, and $h(x) \leq x$ for all $x \in \bR_{+}$.
Then $h(x)=1_{\{x \geq 1\}}$ for all $x \in \bZ_{+}$, and
\begin{eqnarray*}
P(\max_{j \leq m}X_j>u_n)&=&
%P(\sum_{j=1}^{m}1_{\{X_{j} >u_n\}} \geq 1 )=
 E[1_{\{\sum_{j=1}^{m}1_{\{X_{j} >u_n\}} \geq 1\}
}] =E[h(\sum_{j=1}^{m}1_{\{X_{j} >u_n\}})].
\end{eqnarray*}
We omit the details. $\Box$

\vspace{3mm}

\noindent {\bf Proof of Theorem \ref{extremal-th}:} The result
follows from Proposition \ref{calcul-theta} and Proposition
\ref{difference2}, using the fact that: $$P(\max_{j \leq m-1}X_{j}
\leq u_n,X_{m}>u_n)=P(\max_{j \leq m} X_{j}>u_n)-P(\max_{j \leq m-1}
X_{j}>u_n).$$ $\Box$

\vspace{3mm}

\footnotesize{{\em Acknowledgement.} The first author would like to
thank the second author and Laboratoire de math\'ematiques de
l'Universit\'e de Paris-sud (\'Equipe de Probabilit\'es,
Statistique et mod\'elisation) for their warm hospitality.

\normalsize{

}

\end{document}